\input amstex1.tex
\input psfig
\overfullrule=0pt
\pagewidth{4.9in}
\pageheight{7.5in}
\hcorrection{0.53in}
\vcorrection{.02in}
\TagsOnRight
\font\lmathfont=cmsy10 scaled \magstep 1
\long\gdef\ignore#1{}
\define\C{\Bbb C}

\define\R{\Bbb R}

\define\NEarrow{\text{\lmathfont\char'045}}
\define\SEarrow{\text{\lmathfont\char'046}}
\define\NWarrow{\text{\lmathfont\char'055}}
\define\SWarrow{\text{\lmathfont\char'056}}
\define\Z{\Bbb Z}

\def\mapright#1{\smash{\mathop{\longrightarrow}\limits^{#1}}}

\def\mapup#1{\uparrow\rlap{$\vcenter{\hbox{$\scriptstyle#1$}}$}}
\def\mapNW#1{\mathrel{\mathop{\kern0pt \NWarrow}\limits^{\mkern25mu\smash{\lower
6pt\hbox{$\scriptstyle#1$}}}}}
\def\mapSE#1{\mathrel{\mathop{\kern0pt \SEarrow}\limits^{\mkern25mu\smash{\lower
6pt\hbox{$\scriptstyle#1$}}}}}
\def\mapNE#1{\mathrel{\mathop{\kern0pt \NEarrow}\limits^{\mkern-15mu\smash{\lowe
r6pt\hbox{$\scriptstyle#1$}}}}}
\def\mapSW#1{\mathrel{\mathop{\kern0pt \SWarrow}\limits^{\mkern-15mu\smash{\lowe
r
6pt\hbox{$\scriptstyle#1$}}}}}

\redefine\emptyset{/\kern -0.67em{\ssize\bigcirc}}
\def\picture #1 by #2 (#3){\centerline{
  \vbox to #2{
    \hrule width #1 height 0pt depth 0pt
    \vfill
    \special{picture #3} 
    }
    }
  }
\def\scaledpicture #1 by #2 (#3 scaled #4){{
  \dimen0=#1 \dimen1=#2
  \divide\dimen0 by 1000 \multiply\dimen0 by #4
  \divide\dimen1 by 1000 \multiply\dimen1 by #4
  \picture \dimen0 by \dimen1 (#3 scaled #4)}
  }
\def\centercap#1{\write1{\line{#1\string\dotfill\hss\number\pageno}}\centerline{#1}\medskip}

\newif\ifboxfigure      
\boxfigurefalse         

\def\BoxIt#1#2{
	\vbox{\hrule
	\hbox{\vrule\kern#2\vbox{\kern#2#1\kern#2}\kern#2\vrule}
		   \hrule}}

\def\RasterBox #1 #2 #3 #4{


\dimen5=65pt
\divide\dimen5 by 72

\dimen0=#2\dimen5
\divide\dimen0 by 1000
\dimen1=#3\dimen5
\divide\dimen1 by 1000
\dimen2=#3\dimen5
\divide\dimen2 by 1000
\dimen3=#2\dimen5
\divide\dimen3 by 1000

\setbox4=\hbox to #4\dimen0{
 \vbox to #4\dimen1{
 \vss
 \psfig{figure=#1,height=#4\dimen2,width=#4\dimen3}
 }
 \hss
 }
 \ifboxfigure\BoxIt{\box4}{0pt}
 \else\box4
 \fi
 }

\def\insertRaster #1 pixels #2  by #3 scaled #4 {
			\medskip
			 \hbox to \hsize{%

			 \hss
			 \RasterBox {#1} {#2} {#3} {#4}
			 \hss
			 }%
}

\documentstyle{amsppt1}
\nologo

\centerline{\bf Conformal Dynamics Problem List}\medskip
\centerline{\bf Edited by Ben Bielefeld}\medskip
The following list of unsolved problems was given at the Conformal Dynamics 
Conference which was held at SUNY Stony Brook in November 1989.  Problems were contributed
by Ben Bielefeld, Adrien Douady, Curt McMullen, Jack Milnor, Misuhiro
Shishikura, Folkert Tangerman, and Peter Veerman. 
\bigskip
\centerline{\bf \S1. Local connectivity of Julia sets }
\smallskip

Let $f_\lambda(z)=\lambda z + z^2$ where $\lambda=\exp(2\pi i \theta)$.
Call such a polynomial {\it parabolic} if $\theta$ is rational.  If $\theta$
is irrational, and $f_\lambda$ is analytically conjugate to a rotation
near 0 we say $f_\lambda$ is a {\it Siegel polynomial}.  Otherwise we call 
$f_\lambda$ a {\it Cremer polynomial}.  Douady and Sullivan [{\bf Sul}] 
have shown that the Julia set of a
Cremer polynomial is never locally connected. In the generic
case, Douady (unpublished) has described specific examples of external rays which do
not land, but rather have an entire continuum of limit points in the
Julia set.

{\bf Question 1.} Is there an arc joining 0 to $-\lambda$, in the Julia set
of a Cremer polynomial? ( $-\lambda$  is the preimage of the fixed
point.)

{\bf Question 2.} Give a plausible topological model for the Julia set of a 
Cremer polynomial.

{\bf Question 3.} Make a good computer picture of the Julia set
of some Cremer polynomial.

{\bf Question 4.} Are there any rays landing at 0 for a Cremer polynomial?

{\bf Question 5.} For which Cremer polynomials is the critical point accessible?

{\bf Question 6.} If we remove the fixed point from the Julia set of a Cremer
polynomial, how many connected components are there in the
resulting set $J(f_\lambda)-\{0\}\,$, ie., is the number of components
countably infinite?

{\bf Question 7.}  Are there any Siegel polynomials whose Siegel disk has a
boundary which is not a Jordan curve?
\smallskip

Let $P_c(z)=z^2+c$ where $c$ ranges over values for which the Julia set
of $P_c$ is connected.

{\bf Question 8.} For which $c$ is the Julia set of $P_c$ locally connected?
[Reportedly
Yoccoz has recently proved local connectivity except at points on boundaries
of hyperbolic components and infinitely renormalizable points.]

{\bf Question 9.}  Is the Julia set of $P_c$ locally connected when c is real?

{\bf Question 10.}  If $c$ is the Feigenbaum point, is the Julia set of $P_c$ locally connected?

{\bf Question 11.} If one could show that topological conjugacy of $P_{c_1}$
and $P_{c_2}$ implies that $P_{c_1}$ is quasiconformally conjugate to
$P_{c_2}$, would this imply local connectivity of the quadratic
connectedness locus (ie., the Mandelbrot set, consisting of all
 $c$ for which the Julia set of $P_c$ is
connected)?
\smallskip

Let $\tau = p/q$, and let $M_\tau$ be the limb of the connectedness locus
with interior angle $p/q$.

{\bf Question 12.}  Is the diameter of $M_\tau$ less than $K/q^2$
for some constant $K$ independent of $\tau$?  If not, is it
at least less than  $K \log(q)/q^2$? (Remark: the still unpublished
Yoccoz inequality implies that the diameter is bounded by a constant
over  $q$ . Compare [{\bf P}].)

{\bf Question 13.}  If $P_c$ is nonrecurrent does this imply the existence
of a conformal metric (a metric of the form $\rho(z) |dz| $
with integrable singularities) for which $P_c$
is expanding on its Julia set?
[Yoccoz has proved local connectivity in the nonrecurrent
case.]
\bigskip

References:

[{\bf D}] A. Douady, Disques de Siegel et anneaux de Herman, S\'em. Bourbaki
n$^o$ 677, 1986-87.

\ref\key [{\bf DH}1] 
  \by A. Douady and J.H. Hubbard 
  \paper Syst\`emes Dynamiques Holomorphes I,II: It\'eration des Polyn\^omes
Complexes
  \jour Publ. Math. Orsay
  \vol 84.02 and 85.04
\endref

[{\bf G}] E. Ghys, Transformations holomorphes au voisinage d'une courbe de
Jordan, CRAS Paris 298 (1984) 385-388.

\ref\key [{\bf P}]
  \by C. Petersen
  \paper Yoccoz theorem and inequality
  \jour Aarhus Univ.
  \yr in preparation
\endref

\ref\key[{\bf Sul}]
  \by D. Sullivan.
  \paper Conformal dynamical systems
  \jour  Geometric Dynamics,
725--752, Springer-Verlag Lecture
Notes No. 1007, 1983
\endref

[{\bf Y}] J.-C. Yoccoz, Lin\'earisation des germes de diff\'eomorphismes
holomorphes de $({\bold C},0)\,$, CRAS Paris 306 (1988) 55-58.
\bigskip

\centerline{\bf \S2. Quasiconformal Surgery (Douady, Bielefeld, Shishikura)}
\medskip

It is possible to investigate rational functions using the technique of
quasiconformal surgery as developed in [{\bf DH}2], [{\bf BD}] and [{\bf S}].  
There are
various methods of gluing together polynomials via quasiconformal surgery to 
make new polynomials or rational functions.  The idea of 
quasiconformal surgery is to cut and paste the dynamical spaces for two
 polynomials so as to end
up with a branched map whose dynamics combines the dynamics of the 
two polynomials.  One then tries to find a conformal structure that is
 preserved under this branched map of the sphere to itself, so that using
the Ahlfors-Bers theorem the map is conjugate to a rational function.
There are several topological surgeries which experimentally seem to exist,
but
for which no one has yet been able to find a preserved complex structure.

The first such kind
of topological surgery is {\bf mating} of two monic polynomials
with the same degree. (Compare [{\bf TL}].) The first step is to think of each
polynomial as a map on a closed disk
by thinking of infinity as a circle worth of points, one point for each angular
direction.  The obvious extension of the polynomial at the circle at infinity 
is 
$\theta\mapsto d\theta$ where $d$ is the degree of the polynomial.  Now glue
two such polynomials together at the circles at infinity by mapping the $\theta$
of the first polynomial to $-\theta$ in the second. Finally, we must shrink
each of the external rays for the two polynomials to a single point.
The result should be conjugate to a rational map of degree $d$.
(Surprisingly this construction sometimes seems to make sense even when
the filled Julia sets for both polynomials have vacuous interior.)

For instance we can take the rabbit to be the first polynomial, that is 
$z^2 + c$ where the critical point is periodic of period 3 
($c\sim -.122561 + .744862i$).  The Julia set appears in the following
picture.

\insertRaster 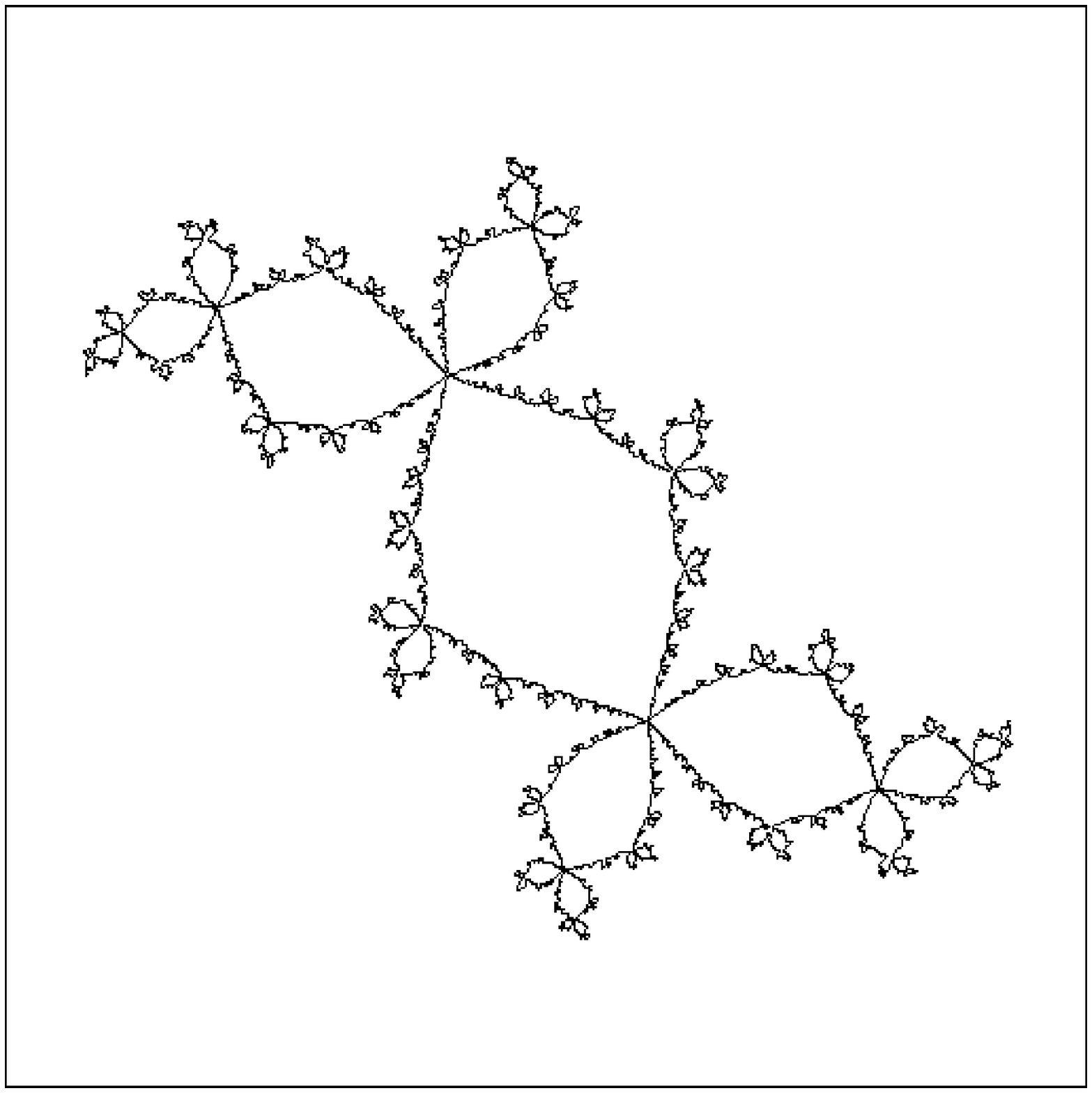 pixels 512 by 512 scaled 250
\centercap{ The rabbit}

Then for the second polynomial we could take the basilica, that is $z^2 - 1$
(it is named after the Basilica San Marco in Venice.  One can see the
basilica on top and its reflection in the water below).
The Julia set for the basilica apears in the following figure.

\boxfiguretrue
\insertRaster 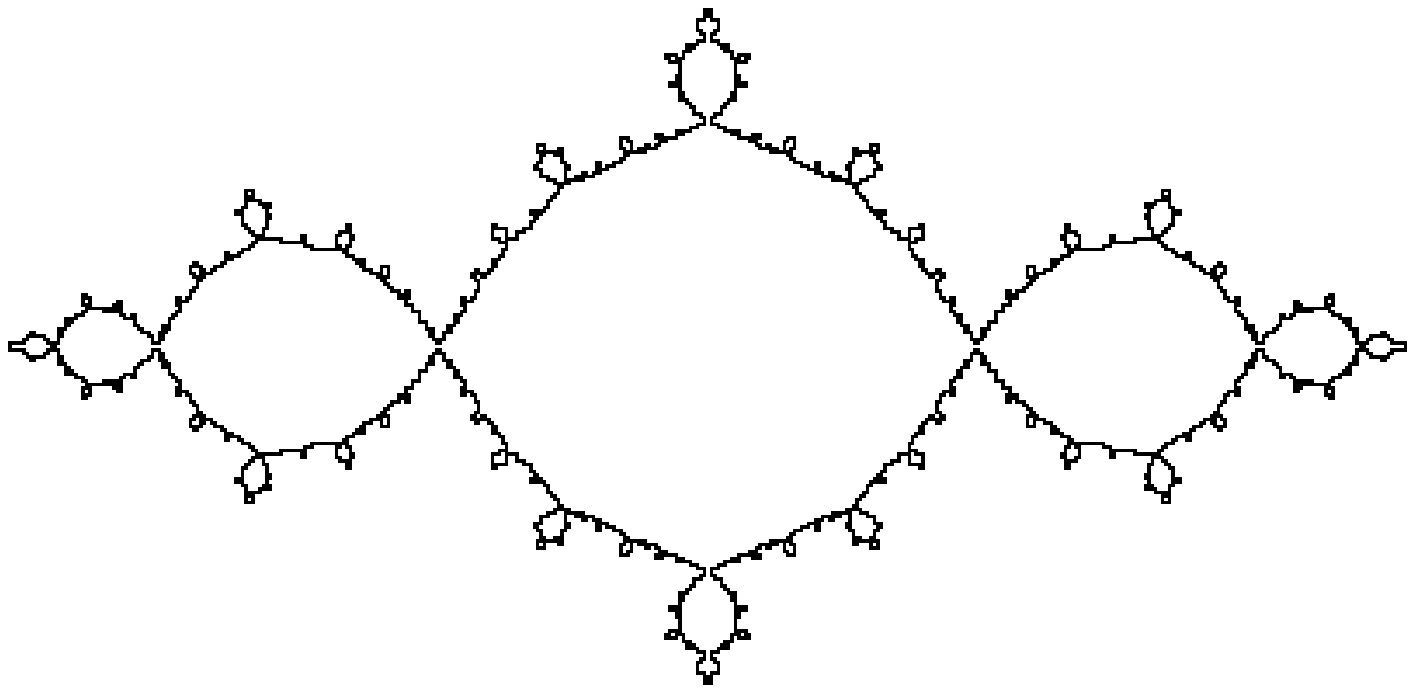 pixels 500 by 500 scaled 250
\centercap{ The basilica}

Next we show the basilica inside-out ($\frac{z^2}{z^2-1}$) 
which is what we will glue to the
rabbit.

\insertRaster 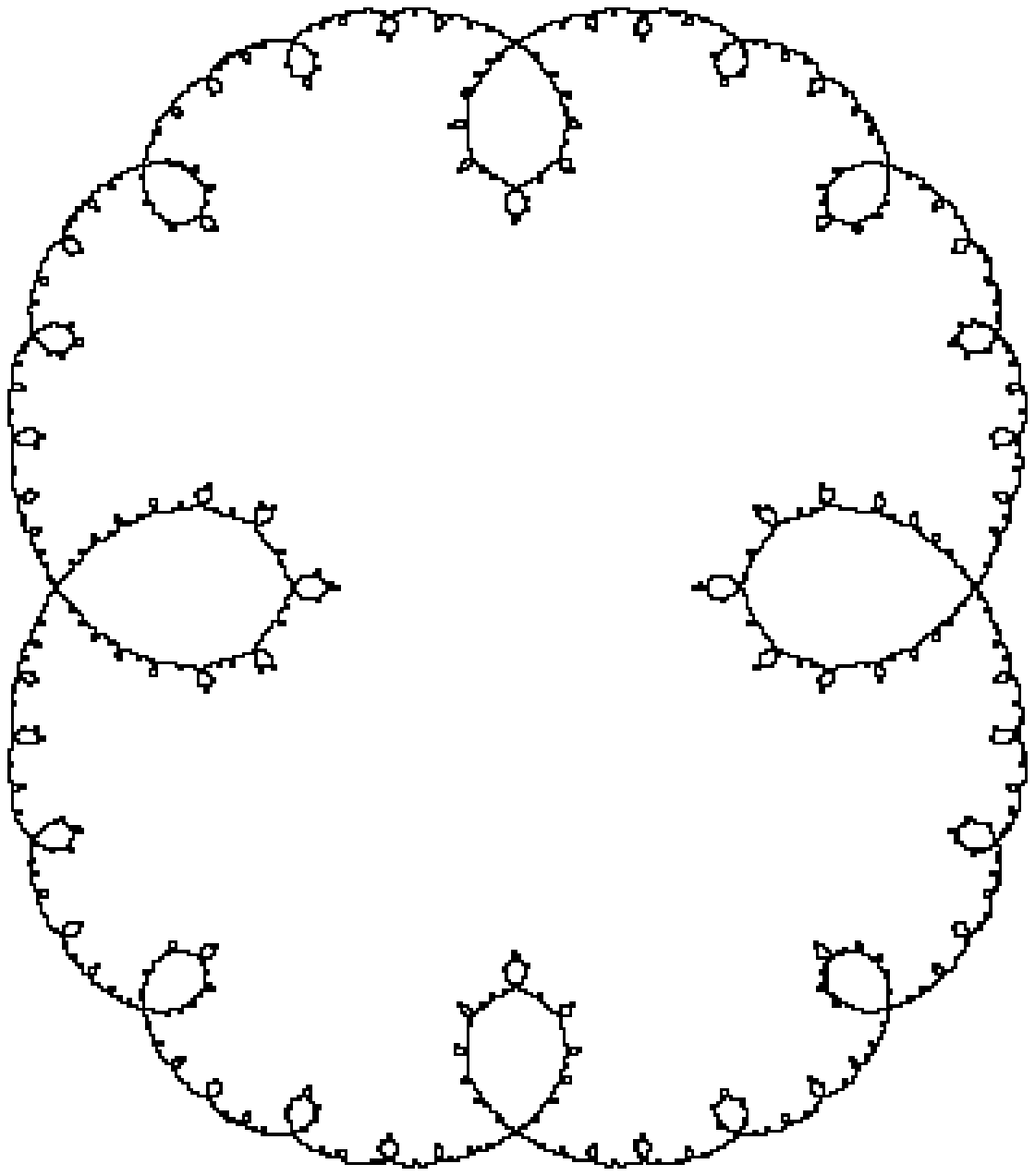 pixels 500 by 500 scaled 250
\centercap{ The inside-out basilica}

And finally we have the Julia set for the mating ($\frac{z^2+c}{z^2-1}$
where $c=\frac{1+\sqrt{-3}}{2}$).

\insertRaster 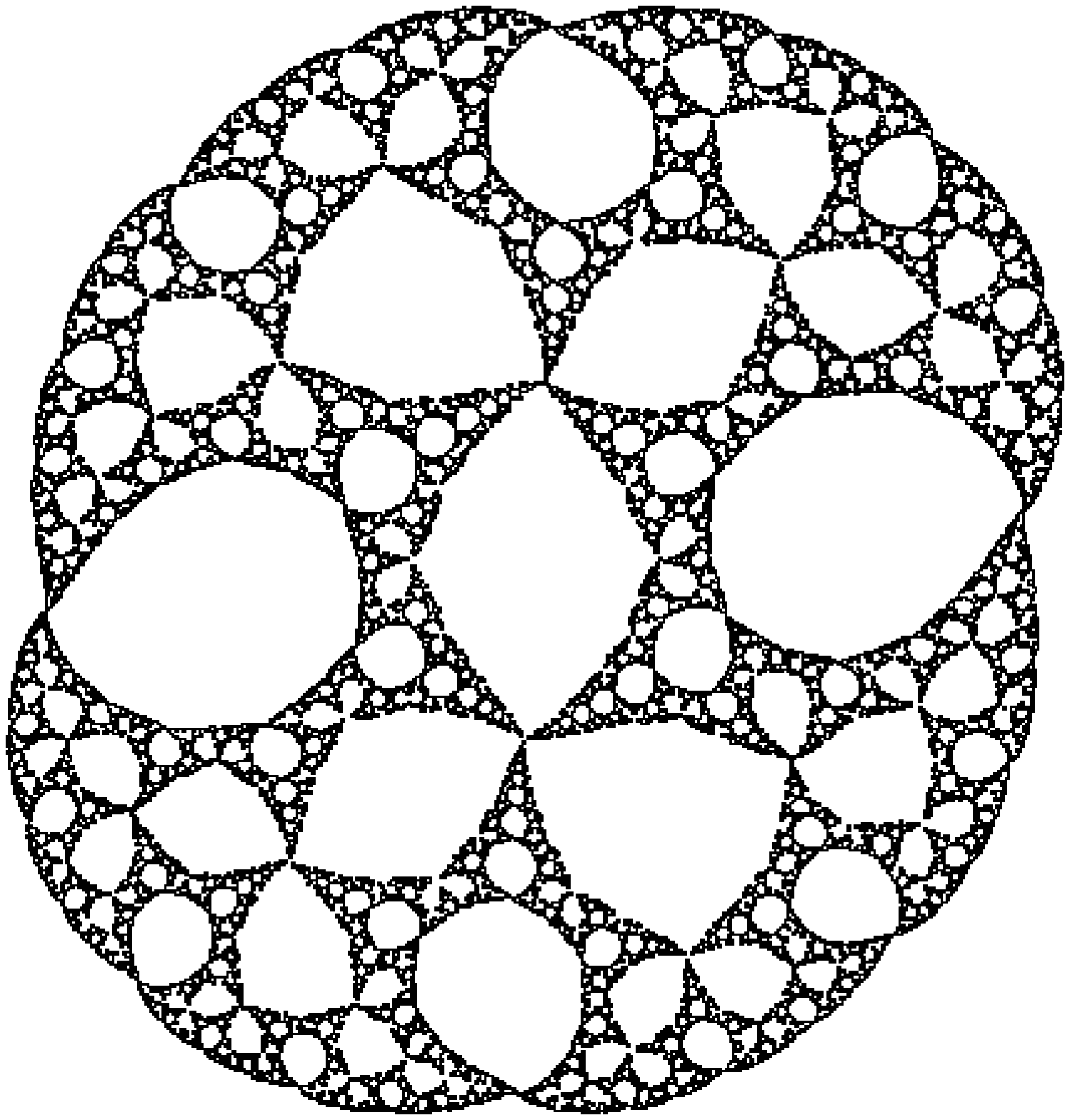 pixels 500 by 500 scaled 250
\centercap{ The basilica mated with the rabbit}
\boxfigurefalse
\eject
{\bf Question 1.}  Which matings correspond to rational functions?  There are
some known obstructions. For example, Tan Lei has shown that matings
between quadratic polynomials can exist only if they do not belong
to complex conjugate limbs of the Mandelbrot set.

{\bf Question 2.} Can matings be constructed with quasiconformal\break
surgery?

{\bf Question 3.}  If one polynomial is held fixed and the other is varied
 continuously, does the resulting rational function vary continuously?
 Is mating a continuous function of two variables?
\smallskip

The second type of topological surgery is {\bf tuning}.
First take a polynomial $P_1$
with a periodic critical point  $\omega$ of period $k$, and assume that no 
other critical points are in the entire basin of this superattractive cycle.   
Let $P_2$ be a polynomial with one critical point whose degree is
the same as the degree of $\omega$.  We also assume that
the Julia sets of $P_1$  and $P_2$ are connected.
We give two descriptions.
For the first description we assume
the closure $\bar B$ of the immediate basin of $\omega$ is homeomorphic
to the closed unit disk $\bar D$,  
and that the Julia set for $P_2$ is locally connected.
Now,
$P_1^k$ maps $\bar B$ to itself by a map
which is conjugate to the map $z\mapsto z^d$ of $\bar D$, where $d$
is the degree of the critical point.
(In fact, if $d>2$, then there are $d-1$ possible choices for the conjugating
homeomorphism, and we must choose one of them.)
Intuitively the idea is now the following. Replace the basin $B$ by a copy
of the dynamical plane for $P_2$, gluing the ``circle at infinity''
for this plane onto the boundary of $B$ so that external angles for
$P_2$ correspond to internal angles in $\bar B$. Now shrink each
external ray for $P_2$ to a point. Also, make an analogous modification
at each pre-image of $B$.  The map from the modified $B$ to 
its image will be given by $P_2$, and the
map on all other inverse images of the modified $B$ will be the identity.
The result,$P_3$, called $P_1$ tuned with
$P_2$ at $\omega$, should be conjugate to a polynomial
having the same degree as $P_1$.  Conversely $P_2$ is said to be
obtained from $P_3$ by renormalization.

In the case of quadratic polynomials, the tunings can be made also
in the case when $P_2$ is not locally connected.  Also it may be
that the case where there is a critical point on the boundary
of the basin is different.

As an example we can take $P_1$ to be the rabbit polynomial. 
Then we can take $P_2(z)=z^2 -2$ which has the closed segment from
-2 to 2 as its Julia set.  The following figure shows the resulting
quadratic Julia set tuning the rabbit with the segment ($z^2+c $ where 
$c\sim -.101096 +.956287i$).

\insertRaster 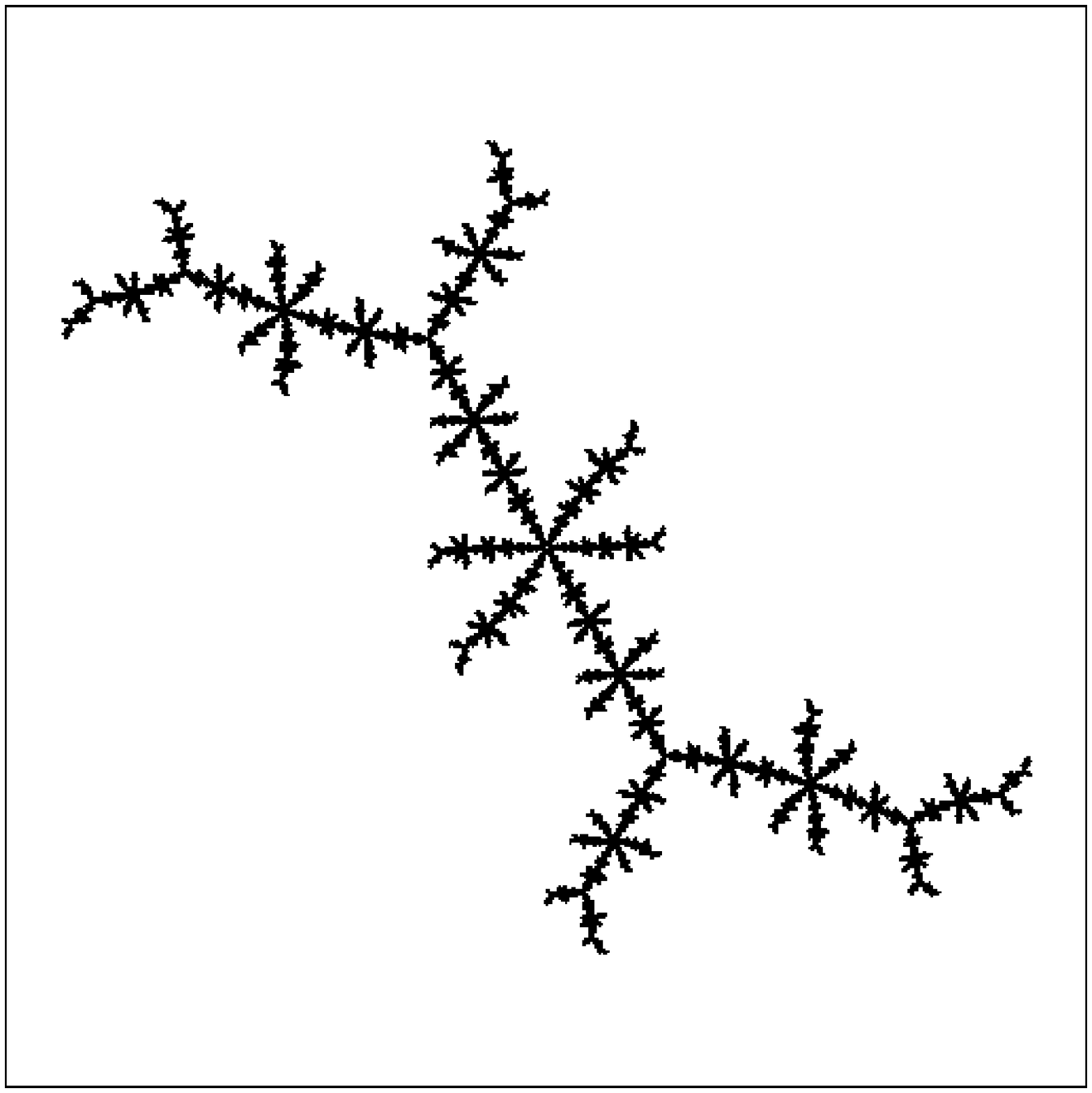 pixels 512 by 512 scaled 250
\centercap{ The rabbit tuned with the segment}

In the picture we see each ear of the rabbit replaced with a segment.

{\bf Question 4.} Does the tuning construction always give a result which is
conjugate to a polynomial?  This is true when
$P_1$ and $P_2$ are quadratic.

{\bf Question 5.} Can tunings be constructed with quasiconformal sur\-gery?  

{\bf Question 6.} Does the resulting polynomial vary continuously with
$P_2$? This is true when $P_1$ and $P_2$ are quadratic [{\bf DH}2].

{\bf Question 7.} Does the resulting tuning vary continuously with $P_1$? (here 
we consider only $P_1$ with a superstable orbit.) 

{\bf Question 8.} Given a sequence $P_{1,k}$ of tending to some limit, do
the tunings of $P_{1,k}$ with $P_2$ tend to a limit which is independent
of $P_2$ ?
\smallskip

The third kind of surgery is {\bf intertwining surgery}.

Let $P_1$ be a monic polynomial with connected Julia set
having a repelling fixed point $x_0$ which has ray landing on it with rotation
number $p/q$.  Look at a cycle of $q$ rays which are the forward images of
the first.
Cut along these rays and we get $q$ disjoint 
wedges.  Now let $P_2$ be a monic polynomial with  a ray of the same
rotation number landing on a repelling periodic point of some period dividing
$q$ (such as 1 or q).
Slit this dynamical plane along the same rays making holes for the wedges.
Fill the holes in by the corresponding wedges above making a new sphere.
The new map will given by $P_1$ and $P_2$ except on a neighborhood of
the inverse images of the cut rays where it will have to be adjusted
to make it continuous.  This construction should be possible to do
quasiconformally using the methods in [{\bf BD}] together with Shishikura's
new (unpublished) method of presurgery in the case where the rays in the
$P_2$ space land at a repelling orbit.  This construction doesn't seem to work
when the rays
land at a parabolic orbit.  

For instance we can take $P_1(z)=z^2$ and $P_2(z)=z^2 - 2$.  The Julia
set for $P_1$ is the unit circle with repelling fixed point at 1 and the
ray at angle 0 lands on it with rotation number 0.  The
Julia set for $P_2$ is the closed segment from -2 to 2 with repelling
fixed point 2 and the
ray at angle 0 lands on it with rotation number 0.
We cut along the 0 ray in both cases.  Opening the
cut in the first dynamical space gives us one wedge.  The space created
by opening the cut in the second space is the hole into which we put
the wedge.  The resulting cubic Julia set is shown in the following
picture (the polynomial is $z^3 + az$ where $a\sim 2.55799i$).

\insertRaster 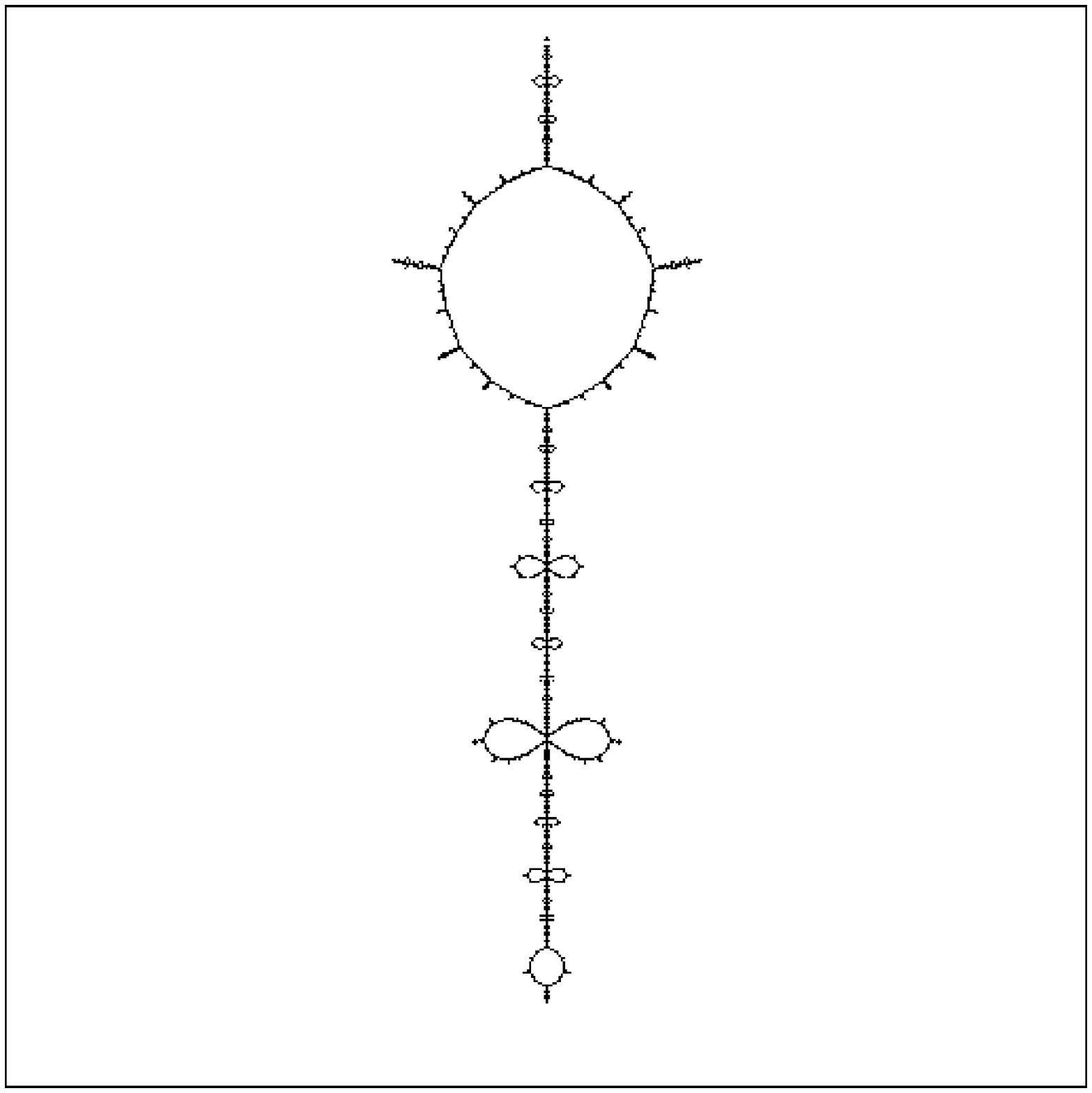 pixels 512 by 512 scaled 250
\centercap{A circle intertwined with a segment}

We see in the picture the circle and the segment, and at the inverse
image of the fixed point on the segment we see another circle. At the
other inverse of the fixed point on the circle we see a segment attached.
All the other decorations come from taking various inverses of the
main circle and segment.

As a second example we can intertwine the basilica with itself. The
ray $1/3$ lands at a fixed point and has rotation number $1/2$.  The
following is the Julia set for the basilica intertwined with itself
(the polynomial here is $z^3 -\frac{3}{4}z + \frac{\sqrt{-7}}{4}$).

\insertRaster 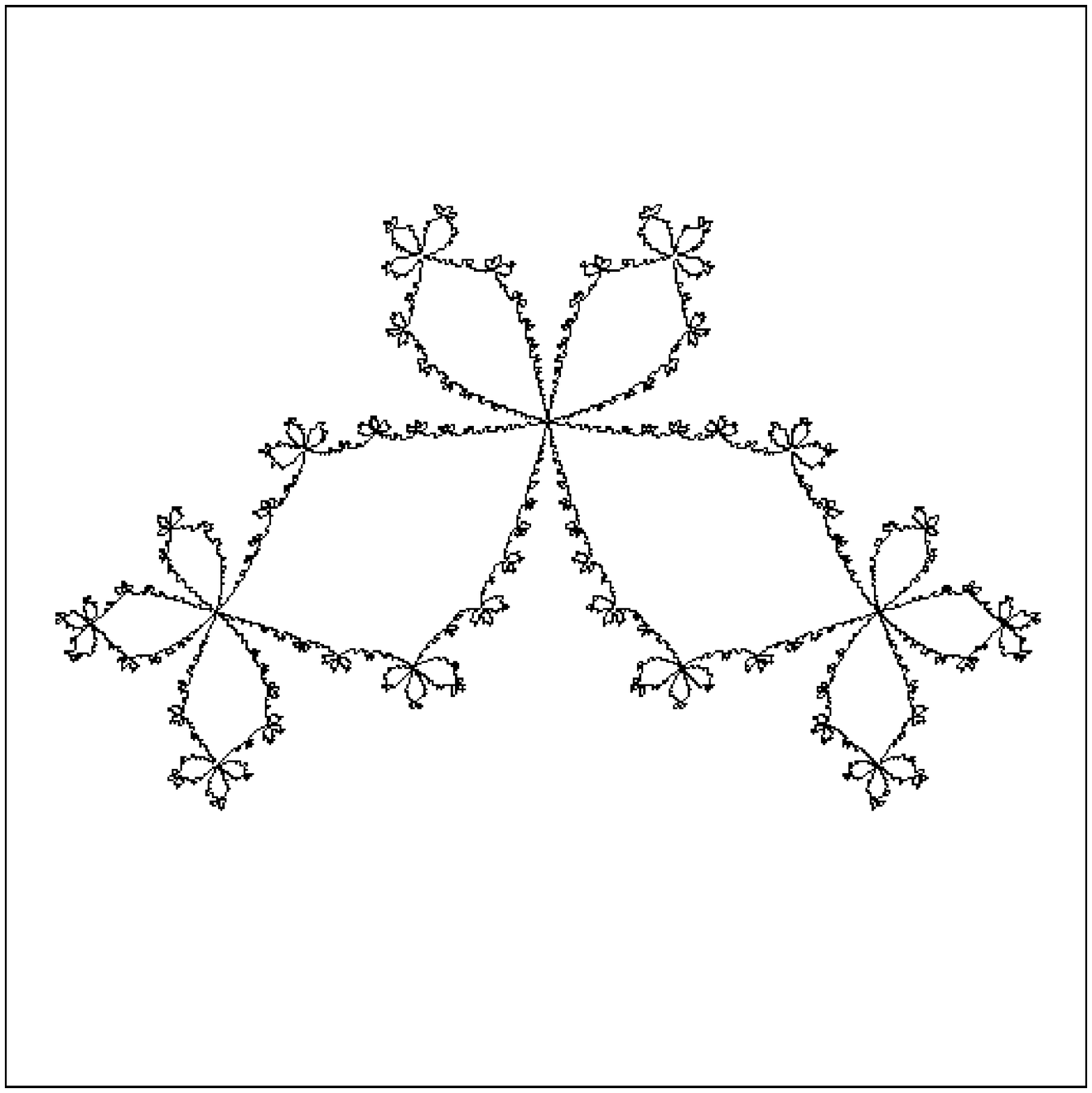 pixels 512 by 512 scaled 250
\centercap{A basilica intertwined with itself}

{\bf Question 9.} When does an intertwining construction give something
which is conjugate to a polynomial?

{\bf Question 10.}  Can intertwinings be constructed with quasiconformal surgery?

{\bf Question 11.}  Does the resulting polynomial vary continuously in $P_2$?
\smallskip

Here is a different kind of continuity question.
Consider the space of all monic polynomials
$$ z \mapsto z^n + a_{n-1} z^{n-1} + \cdots + a_1 z	$$
with  $|a_1| \geq 1\,$, so that there is an
unattractive fixed point at the origin.
Here we do not require that the Julia set be connected.
If at least one external
ray lands at the origin, then there is a well defined ``rotation number" of
these external rays under the map  $\theta \mapsto n \theta$\break (mod 1) .

{\bf Question 12.}  Does this rotation number
extend uniquely to a continuous map from our space of polynomials to  $\R/\Z$ ?
(When  $a_1 = e^{2\pi i\theta}$, this rotation number map must take the value
$\theta$.)

References:

\ref\key [{\bf BD}]
  \by B. Branner and A. Douady 
  \paper Surgery on Complex Polynomials
  \jour Proc. Symp. of Dynamical Systems Mexico
  \yr 1986
\endref

\ref\key [{\bf DH}2] 
  \by A. Douady and J.H. Hubbard 
  \paper On the Dynamics of Polynomial Like Mappings
  \jour Ann. Sc. E.N.S., 4\`eme S\'eries
  \vol 18
  \yr 1966
\endref
\medskip

\ref\key[{\bf S}] 
  \by M. Shishikura 
  \paper On the Quasiconformal Surgery of Rational Functions
  \jour Ann. Sc. E.N.S., 4\`eme S\'eries
  \vol 20
  \yr 1987
\endref
\medskip

\ref\key[{\bf STL}] 
  \by M. Shishikura and Tan Lei
  \paper A Family of Cubic Rational Maps and Matings of Cubic Polynomials
  \jour preprint of Max-Plank-Institute, Bonn
  \vol 50
  \yr 1988
\endref
\medskip

\ref\key[{\bf TL}1]
  \by Tan Lei
  \paper Accouplements des polyn\^omes quadratiques complexes
  \jour CRAS Paris
  \yr 1986
  \pages 635-638
\endref
\medskip

\ref\key[{\bf TL}2]
  \by Tan Lei
  \paper Accouplements des polyn\^omes complexes
  \jour Th\`ese, Orsay
  \yr 1987
\endref
\medskip

\ref\key[{\bf W}]
  \by B. Wittner
  \paper On the Bifurcation Loci of Rational Maps of Degree Two
  \jour Ph.D thesis, Cornell Univ., Ithaca N.Y.
  \yr 1986
\endref
\medskip

\bigskip

\centerline{\bf \S3. Thurston's algorithm for real functions}
\centerline{\bf (Bielefeld,
Tangerman,Veerman,Milnor)}\smallskip

Work in the space
$P = P(n, \mu)$  of piecewise monotone maps  $f$  of the interval $I=[0,1]$ which
have  $n$  laps, and which map the boundary $\{0,1\}$ into itself by some specified
map  $\mu$. Starting with any map  $f_0$  in  $P$ , let  $p_0$ be the unique degree
$n$  polynomial in  $P$  which has all critical points in [0,1] and which has the
same critical values, encountered in the same order, as does  $f_0$. Then there
is a unique homeomorphism  $h_0$  of  [0,1]  which takes the critical points
of  $f_0$  to the critical points of  $p_0$  and which satisfies
$p_0 \circ h_0 = f_0 $.  It follows that the map
  $f_1 = h_0 \circ p_0  = h_0 \circ f_0 \circ h_0^{-1}$  is topologically
conjugate to  $f_0$ .  Now continue inductively, constructing maps
  $f_{k+1} = h_k \circ f_k \circ h_k^{-1} $.
$$
\spreadmatrixlines{5 pt}
\matrix
I & \mapright {f_2} & I \\
\mapup {h_1} & \mapSE {p_1}& \mapup {h_1}\\
I & \mapright {f_1} & I \\
\mapup {h_0} & \mapSE {p_0}& \mapup {h_0}\\
I & \mapright {f_0} & I 
\endmatrix
$$
If  $f_0$  is post-critically finite, and if there is no Thurston obstruction,
then it follows from Thurston's argument that the resulting sequence
of polynomials  $p_k$  converges to a polynomial  $p_\infty$  with the same
kneading sequence as  $f_0$ .  Computer experiments suggest that in many cases
the sequence of maps  $f_k$  converges to this same limit. Furthermore, this
seems to be true even if  $f_0$  is not post-critically finite.

{\bf Question 1.} Formulate and prove some precise statement in this direction.

Now we can consider the same problem except that instead of lifting by
polynomials we lift by some other family in $P$.  For instance computer
experiments suggest that we can choose the $p_k$'s from a family of the form
$\;	p(x) = -|x|^\alpha + c	\;$
where $\alpha > 1\,$.
(To be more precise, we must first change coordinates with an affine
map so that the boundary {0,1} maps to 0.
This yields the family of maps $\;x\mapsto k-k|2x-1|^\alpha\;$ from
the unit interval to itself where $0 < k \leq 1$.)

{\bf Question 2.} Given $f_0$ with a preperiodic or periodic kneading sequence does Thurston's
algorithm converge for the specific family above.  There is a proof
only for $\alpha$ an even integer, which is the polynomial case.

{\bf Question 3.} Give a general property for the lifting family which will
guarantee convergence
of Thurston's algorithm for preperiodic or periodic kneading sequences.

{\bf Question 4.} Give a general property for the lifting family which will
guarantee convergence
of Thurston's algorithm for arbitrary kneading sequences.\bigskip
\vfill\eject
References:

[{\bf MT}] J. Milnor and W. Thurston, Iterated maps of the interval,
pp. 465-563 of ``Dynamical Systems (Maryland 1986-87)'', edit. J.C.Alexan\-der,
Lect. Notes Math. {\bf 1342}, Springer 1988 (cf. \S13.4).
\medskip

\ref\key[{\bf DH}3] 
  \by A. Douady and J.H. Hubbard 
  \paper A Proof of Thurston Topological Characterization of Rational Functions
  \jour Institute Mittag-Leffler Preprint
  \yr  1984
\endref

\bigskip

\def\[{$\,}
\def\]{\,$}
\centerline{\bf \S4. Stable regions for complex H\'enon type maps (Milnor)}
\medskip

Let \[f\] be a polynomial diffeomorphism of \[\C^2\] with Jacobian
determinant \[\delta\]. Suppose that the set \[K^+(f)\] of points with
bounded forward orbit, has a non-empty
interior. Let \[U\] be some connected component of this interior.
According to Montel, the set of iterates of \[f\] restricted to \[U\]
possesses a convergent subsequence, which converges say to \[g:U\to\C^2\].
Evidently the rank \[r\] of \[g\] is zero or one if \[|\delta|<1\],
and is two if \[|\delta|=1\].\smallskip

{\bf Question 1.} Can \[U\] be a wandering component? If so, we must have
\[|\delta|<1\]. Can the rank of $g$ be either zero or one? Can \[U\] be
either bounded, or unbounded of finite volume, or of infinite volume?
\smallskip

If \[U\] is not a wandering component, then it is strictly periodic
under \[f\], and, after replacing \[f\] by some finite iterate, we
may assume that \[f(U)=U\]. Note then that \[f\] commutes with \[g\].
One useful family of examples is provided by the H\'enon
maps. Given any two non-zero complex numbers \[\lambda\] and \[\mu\], there
is an essentially unique (quadratic) H\'enon map \[H_{\lambda,\mu}\] which has a fixed
point with eigenvalues \[\lambda\] and \[\mu\].
\smallskip

{\bf Rank Zero Case.} If \[g(U)={\bold x}_0\in U\], then \[{\bold x}_0\] is an
attracting fixed point. Examples are provided by the H\'enon maps
\[H_{\lambda,\mu}\] where the eigenvalues \[\lambda\] and \[\mu\]
can be any two numbers in the punctured open
disk \[D-0\]. If \[{\bold x}_0\in\partial U\],
then it is conjectured that one eigenvalue must be equal to 1. Evidently
the other eigenvalue must be in \[D-0\]. Here \[H_{1,\mu}\] provides
an example.\smallskip

{\bf Rank One Case.} If \[g\] has a fixed point in \[U\], then \[g\] must be
a retraction of \[U\] onto a Siegel disk. There are examples of the
form \[H_{\lambda,\mu}\] with \[1=|\lambda|>|\mu|\]. (Compare Zehnder.)
Can \[g\] be a retraction onto a Herman ring
(or onto a punctured Siegel disk)?\smallskip

{\bf Rank Two Case.} If \[g\] has a fixed point, then \[U\] is a ``Siegel bi-disk".
There are examples of the form \[H_{\lambda,\mu}\] with \[|\lambda|=
|\mu|=1\]. (Again see Zehnder.)
Can \[U\] also be the product of a
Herman ring with a Herman ring, or the product of a Herman ring
with a Siegel disk?
\bigskip

References: 

[{\bf Z}] E. Zehnder, A simple proof of a generalization of a theorem by C. L.
Siegel, in ``Geometry and Topology III'', ed. do Carmo and Palis, Lecture
Notes Math. {\bf 597}, Springer 1977. 

[{\bf FM}] Friedland and J. Milnor, Dynamical properties of plane polynomial
  automorphisms, Erg. Th. \& Dy. Sy. {\bf 9} (1989), 67-99.

[{\bf HO}] J. Hubbard and R. Oberste-Vorth, H\'enon mappings in the complex domain,
  in preparation.

[{\bf BS}] E. Bedford and J. Smillie, Polynomial diffeomorphisms of $\C^2$:
  currents, equilibrium measure and hyperbolicity, to appear.
\bigskip
\centerline{\bf \S5. Geometrically finite maps and Kleinian groups (McMullen)}
\medskip

\define\arrow{\rightarrow}
\define\Poly{{\rm Poly}}
\define\Rat{{\rm Rat}}
\define\proj{{\Bbb P}}
\define\half{{\Bbb H}}


\centerline{\bf Geometrically finite rational maps}

Let $f(z)$ be a rational map, $C$ its set of critical points,
$P = \cup_1^\infty f^n(C)$ its post-critical set and $J$
its Julia set.
The map $f$ is {\it expanding} if \break $\overline{P} \cap J = \emptyset$.
It is well-known that $f$ is expanding iff some fixed iterate
of $f$ uniformly expands the spherical metric on the Julia set;
these maps are also called {\it hyperbolic} or {\it Axiom A}.

Let the space \Rat$_d$ (respectively \Poly$_d$)
of rational (polynomial) maps of degree $d$ be
equipped with the topology of uniform convergence.
A well-known and fundamental problem is to resolve the following:

{\bf Conjecture}.  The expanding maps form a 
dense subset of \Rat$_d$ and \Poly$_d$.

Cf. [{\bf MSS}] where this is related to the problem of
invariant measurable line fields supported on
the Julia set.  (It is known that the set of expanding
maps is open).  This is not even known in the case of quadratic
polynomials.  The corresponding problems for real maps are also open.

In many ways an expanding rational map is 
well-behaved (cf. [{\bf Sul}]);
it is like a Kleinian group
with compact convex core in $\half^3$.

More generally, let us say a rational map is 
{\it geometrically finite}
if $\overline{P} \cap J$ is a finite set.  Equivalently,
every critical point in the Julia set is pre-periodic. 
(In this case rationally indifferent cycles are allowed).
These maps should be compared to geometrically finite
Kleinian groups.  

For a geometrically finite rational map $f$:

{\bf Problem 1}.
Show that either the Julia set $J$ is the whole sphere
and the action of $f$ on $J$ is ergodic,
or the Hausdorff dimension $\delta$ of $J$ is less than 2.
In the latter case, what can be said about the 
$\delta$-dimensional measure of $J$ and the dynamics with
respect to this measure class?

{\bf Problem 2}.
Show every component of $J$ is locally connected.

{\bf Problem 3}.
Develop for $f$ an analogue of the Haken 
decomposition for 3-manifolds.
For example, if $J$ is disconnected, can $f$ be constructed
by surgery from rational maps with connected Julia sets?

{\bf Problem 4}. 
Extend Thurston's combinatorial theory of critical finite
rational maps (those for which $|P| < \infty$)
to all geometrically finite maps.  That is,
describe $f$ up to combinatorial equivalence rel 
$\overline{P}$ by a finite amount of topological data,
and characterize those combinatorial types which arise
as rational maps.

\centerline	{\bf Convex hyperbolic 3-manifolds}

Let $N$ be a complete hyperbolic 3-manifold presented
as the quotient of $\half^3$ by the action of a Kleinian group
$\Gamma$.
The {\it convex core} of $N$ is the quotient of
the convex hull of the limit set.  
All closed geodesics in $N$ are contained in the convex core.

{\bf Question.}  Suppose $\pi_1(N)$ is generated by $n$
elements.  Is there an upper bound $R_n$ 
to the radius of an embedded ball entirely contained
in the convex core of $N$?
(Here $R_n$ should depend only on $n$).

The question has a positive answer when $N$ is a quasifuchsian
group.  By results of Thurston [{\bf Th Ch. 13}]
there is a pleated surface near every point in the convex
core, and this provides an upper bound on the injectivity radius.

The question also has an (easy) positive solution for
hyperbolic 2-manifolds, and we know of no counterexample
for hyperbolic manifolds of any dimension.

\centerline	{\bf Critically finite rational maps on $\proj^n$}

A basic tool aiding the study of critically finite rational 
maps on the Riemann sphere is the Poincar\'e metric on
the complement of the post-critical set $P$ (assuming $|P| > 2$).
This metric is expanded by $f$.
One should be able to apply
the same sort of arguments to critically finite rational maps 
$f : \proj^n \arrow \proj^n$, $n > 1$, 
such that the complement of the
post-critical set is Kobayashi hyperbolic. 

More precisely, say $f$ is {\it critical finite} if
there exist (possibly reducible) hypersurfaces
$V \subset W \subset \proj^n$ such that 
$$
	f : (\proj^n - W) \arrow (\proj^n - V)
$$
is a covering map.  

{\bf Problem.}  Are there nontrivial examples of
critically finite maps with $\proj^n-V$ Kobayashi
hyperbolic?  How do they behave dynamically?  

\bigskip
References:
\medskip

\ref\key[{\bf MSS}]
  \by R. Ma\~ne, P. Sad, and D. Sullivan.
\paper On the dynamics of rational maps
 \jour{\it Ann. Sci. \'Ec. Norm. Sup.} 16 (1983), 193--217
\endref

\ref\key[{\bf Sul}]
  \by D. Sullivan.
\paper Conformal dynamical systems
\jour  Geometric Dynamics,
 725--752, Springer-Verlag Lecture
 Notes No. 1007, 1983
 \endref

\ref\key[{\bf Th}]
  \by W.~P. Thurston.
\paper  Geometry and Topology of Three-Manifolds
\jour Princeton lecture notes, 1979
 \endref


\end